\documentclass{article}
\usepackage{graphicx,amsmath,amsthm,amsfonts} 
\usepackage{todonotes}
\usepackage{hyperref}

\newtheorem{theorem}{Theorem}

\theoremstyle{definition}
\newtheorem{quest}{Question}

\numberwithin{equation}{section}

\newcommand{\Z}{\mathbb{Z}}

\newcommand{\R}{\mathbb{R}}

\newcommand{\Sc}{\mathbb{S}^1}

\newcommand{\Diff}{\mathrm{Diff}}

\newcommand{\bp}{p}
\newcommand{\bq}{q}

\newcommand{\frp}{P}
\newcommand{\frq}{Q}

\title{Denjoy examples of class $C^1$ with affine dynamics outside the invariant Cantor set}
\author{Victor Kleptsyn, Andrés Navas}
\date{}

\begin{document}

\maketitle

\begin{abstract}
    Given any irrational number, we construct a Denjoy example of class $C^1$ with this rotation number, that is exactly affine on each interval of the complement to the exceptional minimal set.
\end{abstract}

\section{Introduction}

Answering a question of Poincaré, Denjoy~\cite{Denjoy} proved in 1932 that 
every $C^2$ circle diffeomorphism with an irrational rotation number is 
minimal, though this is no longer true for $C^1$ diffeomorphisms. Since 
then, circle maps with an irrational rotation number having an invariant 
Cantor set have been known as {\em Denjoy (counter-)examples} (yet Poincaré 
himself knew of their existence in the continuous setting~\cite{Poincare}, 
and Bohl had previously achieved the construction of one of them in the 
$C^1$ category~\cite{Bohl}).

In his seminal work~\cite{Herman}, Herman revisited this 
matter. In particular, in Chapter X, \S 3 therein, 
he proved that for every irrational angle $\alpha$, 
there exists a Denjoy example $f$ with rotation number $\alpha$ 
that is of class $C^{1+\tau}$ for every $0\!<\!\tau\!<\!1$. Let us emphasize 
that, in this construction, the derivative of $f$ is identically equal 
to~$1$ on the invariant Cantor set (and, in particular, at the endpoints 
of the intervals of its complement). 

Despite being a classical topic that has been thoroughly studied, several 
questions remain open. Here we list three of them that, to our taste, are 
among the most significant.

\begin{quest}
    Is Denjoy's theorem still true for a diffeomorphism whose derivative is $\omega$-continuous with respect to the modulus of continuity $\omega(s) = s \lvert\log(s)\rvert$?
\end{quest}

This question has been raised by Katok, among others (see \cite[\S 3.12]{KH}). 
Note that the modulus of continuity $\omega$ above is weaker than Lipschitz but stronger than 
$\tau$-Hölder for every $\tau < 1$. Moreover, an elementary argument shows that no Denjoy example with an $\omega$-continuous derivative can be built with derivative identically 
equal to $1$ on the invariant Cantor set (see \cite[Exercise 4.1.26]{Na-book}).

\begin{quest}
Given an irrational angle $\alpha$, what is the optimal modulus of continuity $\omega$ 
for which $C^{1+\omega}$ circle diffeomorphisms with irrational rotation number $\alpha$ are necessarily minimal?  
\end{quest}

In other words, what is the best modulus of continuity for which Denjoy's theorem holds for maps with a prescribed 
rotation number $\alpha$? This question is implicit in~\cite{Herman}. In particular, in Chapter X, \S 4 therein, Herman 
proved that for extremely Liouville numbers $\alpha$, this optimal modulus is strictly weaker than Lipschitz.

\begin{quest}\label{q:McDuff}
For a $C^1$-Denjoy example $f$, what can be said about the distribution of the lengths of the intervals in the complement of the invariant Cantor set? In concrete terms, if we let $\ell_1$ be the maximum of these lengths and inductively define $\ell_n$ as the maximum among the lengths that are smaller than $\ell_{n-1}$, does the quotient $\ell_{n-1}/\ell_n$ converge to $1$?
\end{quest}

This question was raised by McDuff in~\cite{McDuff}, where she proved, for instance, that the ternary Cantor set cannot 
arise as the invariant set of a $C^1$-Denjoy example (see \cite{Kercheval} and \cite{Portela} for more on this). 

\vspace{0.2cm}

The purpose of this note is to pursue this topic by exhibiting somewhat surprising examples of Denjoy examples. More 
concretely, we construct $C^1$ Denjoy examples that are affine on each interval of the complement of the invariant Cantor set.

\begin{theorem}\label{t:main}
For any irrational $\alpha$, there exists $f \in \Diff^1_+(\mathbb{S}^1)$ with rotation number $\rho(f)=\alpha$ that has an invariant 
Cantor set~$K$ such that, when restricted to any interval in the complement $\mathbb{S}^1 \setminus K$, the map $f$ is affine.
\end{theorem}

We wrote this note with the hope that this result would shed light on the problems listed above, as well as many others concerning critical regularity for circle and interval maps (in this direction, see for instance \cite{DKN}). And actually, by the time this article was completed, using a brilliant modification of the present construction, Maximiliano Escayola gave an example providing a (very unexpected) negative answer to Question~\ref{q:McDuff}. 
We thus urge the reader, after reading this note, to go through to his work \cite{Escayola}.


\section{The construction}

\subsection{A quotient circle}

We will construct the diffeomorphism $f$ simultaneously with its semi-conjugacy to the circle rotation $R_{\alpha}$. 
Namely, we take the orbit $\{ \{n\alpha\} : {n\in\Z} \} \subset \Sc$ of the point $0$ under $R_{\alpha}$, and we replace each 
point $x_n := \{n \alpha \}= n\alpha \,\, (\mathrm{mod } \,\, \Z)\in \Sc$ by an interval of a certain length $l_n>0$. These lengths $l_n$, to be chosen 
and fixed later, will satisfy the finiteness condition 
\begin{equation}\label{eq:finite}
    L:=\sum_{n\in\Z} l_n <\infty.
\end{equation}

The complement of the union of interiors of these intervals, that is, the minimal invariant set $K$, is declared to be of 
zero Lebesgue measure. In this way, we obtain the circle $\R/(L\Z)$ together with the degree-one (semiconjugacy)  
map $p:\R/(L\Z)\to \R/\Z$ defined by
\[
p(t)=\sup\left\{ x\in [0,1) \mid \sum_{n\in\Z, \, \{n\alpha\} \in [0,x)} l_n \le t\right\}.
\]
For each $n$, this map $p$ projects the interval $I_n:=p^{-1}( x_n )$ of length $l_n$ to the point $x_n \in \R/\Z$.

\begin{figure}
    \centering
    \includegraphics[width=0.5\linewidth]{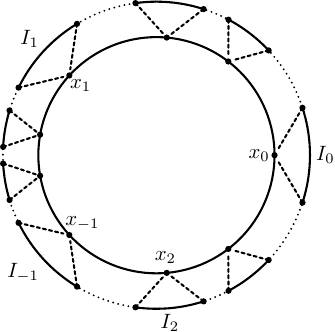}
    \caption{The construction of an affine Denjoy example. 
    The intervals $I_n$ are ``glued'' to the points $x_n=\{n\alpha\}$, $n\in \Z$, and are cyclically permuted in an affine way.}
    \label{fig:Denjoy}
\end{figure}

There exists a unique homeomorphism $f:\R/L\Z \to \R/L\Z$ such that, for each $n \in \mathbb{Z}$, the restriction $f:I_n\to I_{n+1}$ 
is an affine map. The map $p$ is indeed  a semi-conjugacy between $f$ and the rotation $R_{\alpha}$ 
(see Fig.~\ref{fig:Denjoy}).

\subsection{The derivatives and the lengths of the interval}
By construction, the following condition must be satisfied:
\[
\forall n\in\Z, \quad \forall x\in I_n: \qquad \log Df(x)= \log \frac{l_{n+1}}{l_n} = \log l_{n+1} - \log l_n.
\]
In other words, $f$ is a $C^1$ diffeomorphism if and only if these values of the logarithmic derivative 
extend to a continuous function on the circle. Since $\log Df$ is constant inside each interval~$I_n$, 
this is equivalent to that the function 
\[
    \varphi:\{ n \alpha \,\,\,  (\!\!\!\!\!\! \mod \Z) \}_{n\in\Z} \to \R, \quad \varphi( \{ n\alpha \}) = \log l_{n+1} - \log l_n,
\]
extends to a continuous function on the circle $\Sc=\R/\Z$.

We will construct the function $\varphi(x)$ by its Fourier series. To do this, we let 
\begin{equation}\label{eq:varphi-Fourer}
    \varphi(x)= \sum_{k\in\Z} c_k e^{2\pi i k x} 
\end{equation}
for some coefficients $c_k$ to be fixed. These will satisfy several properties, starting with the convergence condition 
\begin{equation}\label{eq:c_k-finite}
\sum_{k\in\Z} |c_k|<\infty.     
\end{equation}
This ensures that the series~\eqref{eq:varphi-Fourer} uniformly converges and thus defines a continuous function on the circle.

\subsection{The Fourier series of the semiconjugacy}

Instead of fixing the series~\eqref{eq:varphi-Fourer} directly, we will define the lengths $l_n$ by 
\[
\log l_n=\psi(n\alpha),
\]
where $\psi(x)$ is given by a (formal) series
\begin{equation}\label{eq:psi-Fourier}
    \psi(x) := \sum_{m\geq 1} d_m \big( \cos (2\pi \frq_m x) -1 \big)
\end{equation}
for a certain sequence of frequencies $(\frq_m)_{m\in \mathbb{N}}$ and a sequence of positive numbers $d_m>0$ 
to be chosen later, in a way that the series~\eqref{eq:psi-Fourier} converges at all points $x_n = n\alpha \,\,  (\mathrm{mod } \,\, \Z)$, where $n\in\Z$.

This choice leads to defining $\varphi$ at the points $x_n$ as 
\begin{multline}\label{eq:phi-from-psi}
    \varphi(x_n)= \log l_{n+1}- \log l_n = \psi(n\alpha+\alpha) - \psi(n\alpha) 
\\
    = \sum_{m\geq 1} d_m \left[\cos (2\pi \frq_m x_n+2\pi \frq_m\alpha) - \cos (2\pi \frq_m x_n)\right].
\end{multline}
Take $\frp_m\in \Z$ such that $|\frq_m \alpha -\frp_m|<\frac{1}{2}$, and let $\alpha_m:=\frq_m \alpha -\frp_m$. 
Then $2\pi \frq_m\alpha$ in the right hand side of~\eqref{eq:phi-from-psi} can be replaced by $2\pi \alpha_m$, leading to
\[
\varphi(x_n)= -\sum_{m\geq 1} 2 d_m \sin (\pi \alpha_m) \cdot 
\sin \left( 2\pi \big( \frq_m x_n + \frac{\alpha_m}{2} \big) \right).
\]
This forces the Fourier series for the function $\varphi$ to be chosen as 
\begin{equation}\label{eq:varphi-sine}
    \varphi(x) = -\sum_{m\geq 1} 2 d_m \sin (\pi \alpha_m) \cdot 
\sin \left( 2\pi \big( \frq_m x +\frac{\alpha_m}{2} \big) \right).
\end{equation}

By rewriting the series~\eqref{eq:varphi-sine} in the form~\eqref{eq:varphi-Fourer}, one 
easily checks that the convergence condition~\eqref{eq:c_k-finite} is satisfied whenever  
\begin{equation}\label{eq:d-convergence}
    \sum_{m \geq 1} d_m \cdot \left|\alpha_m \right| < \infty.
\end{equation}
Obviously, in order to satisfy this condition~\eqref{eq:d-convergence}, 
one can take any positive convergent series $\sum_{m\geq 1} b_m<\infty$ of positive numbers and let
\begin{equation}\label{def:d-eme}
   d_m := \frac{b_m}{|\alpha_m|}.
\end{equation}
This choice ensures that $\varphi$ is a well-defined, continuous function. Note that the series~\eqref{eq:psi-Fourier} 
for $\psi(x)$ converges at $x=0$ (the right-side expression is identically zero). Moreover, recall that the value of  
$\varphi(n\alpha)$ was fixed by~\eqref{eq:phi-from-psi} as the increment from $\psi(n\alpha)$ to $\psi((n+1)\alpha)$. 
Thus, the convergence of the series defining $\varphi$ implies that the series~\eqref{eq:psi-Fourier} converges at the points 
$x_n=n\alpha$ for all $n\in\Z$, and the lengths $l_n=\exp(\psi(n\alpha))$ are therefore well-defined. It now suffices to ensure 
the convergence of the series~\eqref{eq:finite} representing the sum of the interval lengths. This will be achieved in the 
next sections by appropriately choosing the frequencies~$\frq_m$ and the corresponding coefficients~$b_m$. For the sake of the exposition, we will briefly recall some generalities on good approximants in Sec.~\ref{s:contfrac}, consider the model case of $\alpha$ being the golden mean in Sec.~\ref{s:golden}, and complete the construction for a general irrational angle~$\alpha$ in Sec.~\ref{s:general}.


\subsection{Preliminaries: good approximants}\label{s:contfrac}

\begin{figure}[h!]
    \centering
    \includegraphics[width=0.9\linewidth]{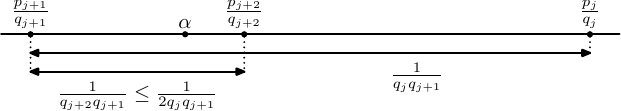}
    \caption{Three consecutive approximantions $\frac{\bp_j}{\bq_j}$ for $\alpha$ and the upper bounds for the distances.}
    \label{fig:cont-frac}
\end{figure}

Let $\frac{\bp_j}{\bq_j}$ be the sequence of good approximations (truncated continued fractions) for the number $\alpha$. 
It is very well-known that $\alpha$ lies between any two consecutive approximants $\frac{\bp_j}{\bq_j}$ and $\frac{\bp_{j+1}}{\bq_{j+1}}$, and that 
\[
\frac{\bp_j}{\bq_j} - \frac{\bp_{j+1}}{\bq_{j+1}}= \frac{(-1)^j}{\bq_j \bq_{j+1}}.
\]
This equality immediately implies an upper bound for the distance to the approximant: 
\begin{equation}\label{eq:c-f-upper}
|\alpha-\frac{\bp_j}{\bq_j}|< \frac{1}{\bq_j \bq_{j+1}}.    
\end{equation}
Moreover, as the next approximant, $\frac{\bp_{j+1}}{\bq_{j+1}}$, is closer to $\alpha$ than $\frac{\bp_j}{\bq_j}$, one has $|\alpha-\frac{\bp_j}{\bq_j}|> \frac{1}{2\bq_j \bq_{j+1}}$. Another way to establish this lower bound is by noticing that $\bq_{j+2}\ge \bq_{j+1}+\bq_{j}\ge 2\bq_{j}$, and thus
\begin{equation}\label{eq:c-f-lower}
\left|\alpha-\frac{\bp_j}{\bq_{j}} \right| > \frac{1}{\bq_j \bq_{j+1}} - \frac{1}{\bq_{j+1}\bq_{j+2}} \ge \frac{1}{2\bq_j\bq_{j+1}}
\end{equation}
see Fig.~\ref{fig:cont-frac}.

Combining the upper and the lower bounds~\eqref{eq:c-f-upper}, \eqref{eq:c-f-lower}, and multiplying by $q_j$, one gets
\begin{equation}\label{eq:interval}
\frac{1}{2\bq_{j+1}} <|\bp_j-\bq_{j} \alpha|< \frac{1}{\bq_{j+1}}.    
\end{equation}

\subsection{A particular case: $\alpha$ is the golden mean}\label{s:golden}

We start with the case where $\alpha=\frac{\sqrt{5}-1}{2}$. In this case, the denominators $q_j$ are the Fibonacci numbers. We will use their exact values as frequencies, setting 
\[
Q_m := q_m, \quad b_m:=\frac{1}{m^2}.
\]
Then, the series $\sum_m b_m$ converges. Taking~$d_m$ to be defined by~\eqref{eq:d-convergence}, we see that the function~$\varphi$ defined by~\eqref{eq:varphi-sine} is well-defined and continuous everywhere on the circle. It thus suffices to verify that the corresponding lengths~$l_n$ form a convergent series. 

\begin{figure}[h!]
    \centering
    \includegraphics[width=0.9\linewidth]{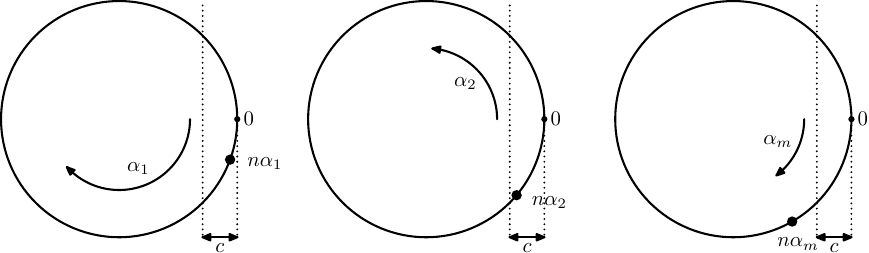}
    \caption{A family of circles equipped with the rotations by $\alpha_m$, and a family of $n$-th iterations of the point~$0$ for a particular value of $n$ (bold points).}
    \label{fig:Rotations}
\end{figure}

Note that each of the summands in the right-hand side expression of~\eqref{eq:psi-Fourier} is nonpositive. Hence, for every $n$ and $m$, one has 
\begin{equation}\label{eq:l-upper-0}
l_n \le \exp(d_m \cdot (\cos (2\pi Q_m \cdot n \alpha) -1)) 
= \exp(d_m \cdot (\cos (2\pi n \alpha_m) -1)).
\end{equation}
As an illustration, one can interpret the angles $2\pi n \alpha_m$, occurring in the right hand side of~\eqref{eq:l-upper-0} in the following way: We consider a countable family of circles, with the $m$-th circle being equipped with the rotation by~$\alpha_m$. These angles become smaller and smaller as $m\to\infty$. The angles appearing in the right-hand hand side of~\eqref{eq:l-upper-0} for a given $n$ correspond to the images of the point $0$ on these circles under $n$ simultaneous rotations; see Fig.~\ref{fig:Rotations}. We will show that at least one of these angles is always sufficiently far from~$0$.

    Namely, note that for each~$m \geq 1$ and $|n|\in [\frac{q_{m+1}}{4},\frac{q_{m+1}}{2}]$, by~\eqref{eq:interval} one has $|\alpha_m n| \in [\frac{1}{8},\frac{1}{2}]$. Hence, 
    \[
        \cos (2\pi n \alpha_m) \le \frac{1}{\sqrt{2}},
    \]
    which implies that  
    \begin{equation}\label{eq:c-lower}
        1-\cos 2\pi n \alpha_m \ge c, \quad \text{ where } \, c:=1-\frac{1}{\sqrt{2}}>0.     
    \end{equation}
    Thus, for such $n$ and $m$, we have the upper estimate 
    \begin{equation}\label{eq:l-upper-golden}
        l_n \le \exp (-c d_m) = \exp \left( - c \frac{1}{m^2 |\alpha_m|} \right) \le \exp \left( - c \, \frac{q_{m+1}}{m^2} \right).   
    \end{equation}
    
    The quotient between consecutive Fibonacci numbers does not exceed~$2$.  Therefore, every number~$n>0$ is contained in at least one semi-open interval $\left[\frac{q_{m+1}}{4},\frac{q_{m+1}}{2}\right[$. We thus have 
    \begin{multline}\label{eq:sum-exp}
        \sum_{n\in \Z} l_n \le 1+ \sum_{m=1}^{\infty}  \, \sum_{|n|\in \left[\frac{q_{m+1}}{4},\frac{q_{m+1}}{2}\right[}  \exp \left(- c \, \frac{q_{m+1}}{m^2} \right) \\ 
        \le 1+ \sum_m \, 2 \cdot \frac{q_{m+1}}{4}  \cdot \exp \left( - c \, \frac{q_{m+1}}{m^2} \right).
    \end{multline}
    The sum in the right hand side of~\eqref{eq:sum-exp} can easily be seen to converge due to the exponential growth of the Fibonacci numbers. This concludes the proof of Theorem~\ref{t:main} in this particular case.

\subsection{The general case: the choice of the frequencies $\frq_m$}\label{s:general}

The argument of Section~\ref{s:golden} can be modified to handle any angle  $\alpha$ of bounded type, at the cost of changing the constant~$c$ in the estimate~\eqref{eq:c-lower}. However, for a general $\alpha$, the quotient of the consecutive denominators $\frac{q_{j+1}}{q_j}$ can be arbitrarily large. Hence, if we were following the same choice $Q_m=q_m$, it might happen that for some $n$ on some of the circles on Figure~\ref{fig:Rotations} after $n$ rotations starting point~$0$ already has returned close to itself, while on the others it even didn't start leaving its neighbourhood.

To overcome this difficulty, we include additional frequencies in the sequence~$(Q_m)$. Namely, for every $j$, consider the sequence of doubled frequencies
\begin{equation}\label{eq:double-j}
\bq_j, 2\bq_j, \dots, 2^{k_j} \bq_j,     
\end{equation}
where 
\[
k_j:=\max \big\{ k \mid 2^k \bq_j < \bq_{j+1} \big\} = \left[ \log_2 \frac{\bq_{j+1}}{\bq_j} \right].
\]
Now, let $(\frq_m)$ be the sequence obtained by concatenating finite sequences~\eqref{eq:double-j}:
\[
(\frq_m) = (\bq_1,\dots, 2^{k_1}\bq_1, \bq_2, \dots, 2^{k_2}\bq_2, \dots, \bq_j, \dots, 2^{k_j}\bq_j, \dots),
\]
and choose 
\begin{equation}\label{eq:b-m}
    b_m = \frac{2^{-(k_j-r)/2} }{j^2} \quad \text{ if }\, \frq_m =2^r \bq_j.
\end{equation}
With this choice, we have 
$$\sum_{m \geq 1} b_m \leq \sum_{k \geq 0} \frac{1}{2^{k/2}} \cdot  \sum_{j \geq 1} \frac{1}{j^2} < \infty.$$
As before, to conclude the proof of Theorem~\ref{t:main}, it suffices to show that, when letting $d_m$ be defined by 
(\ref{def:d-eme}), the sum of the corresponding lengths~$l_n$ converges. 

Again, due to the nonpositivity of the summands in the right-hand side of~\eqref{eq:psi-Fourier}, for every $n$ and $m$ one has 
\begin{equation}\label{eq:l-upper}
l_n \le 
\exp(d_m \cdot (\cos (2\pi n \alpha_m) -1)).
\end{equation}
Now, for $\frq_m=2^r \bq_j$ (with $0 \leq r \leq k_j$), one gets from~\eqref{eq:interval} that 
\begin{equation}\label{eq:alpha-bounds}
|\alpha_m| = |2^r(\bq_j \alpha - \bp_j)| \in \left[ \frac{2^{r-1}}{\bq_{j+1}},\frac{2^{r}}{\bq_{j+1}} \right].    
\end{equation}
For $|n|\in [\frac{\bq_{j+1}}{2^{r+2}},\frac{\bq_{j+1}}{2^{r+1}}]$, one has $|n\alpha_m|\in [\frac{1}{8},\frac{1}{2}]$. In the same way as before (see Eq.~\eqref{eq:l-upper-golden}), substituting 
\[
\cos (2\pi n \alpha_m) \le \frac{1}{\sqrt{2}} 
\]
by the upper bound~\eqref{eq:l-upper}, one gets
\begin{equation}\label{eq:estimate-ln}
l_n \le \exp\left(-\frac{c \, b_m}{|\alpha_m|} \right), \quad \mbox{where } \, c =  1-\frac{1}{\sqrt{2}}> 0.
\end{equation}
Combining the choice~\eqref{eq:b-m} for $b_m$ with the lower bound from \eqref{eq:alpha-bounds} for $\alpha_m$ 
and the definition of $k_j$, we get
\[
\frac{b_m}{|\alpha_m|} \ge \frac{2^{-(k_j-r)/2} }{j^2} \cdot \frac{q_{j+1}}{2^{r}} \ge 2^{(k_j-r)/2} \frac{q_j}{j^2}.
\]
By introducing this inequality into (\ref{eq:estimate-ln}), we obtain the upper bound
\[
    l_n \le \exp\left(- \frac{c \, q_j}{j^2} \cdot 2^{(k_j-r)/2}  \right) \quad \mbox{for } \, |n|\in \left[ \frac{\bq_{j+1}}{2^{r+2}},\frac{\bq_{j+1}}{2^{r+1}} \right[.
\]
There are at most $2\frac{\bq_{j+1}}{2^{r+2}} \le 2^{k_j-r}\bq_{j}$ such indices $n\in \Z$, hence the total contribution of such intervals does not exceed 
\begin{equation}\label{eq:l-q-r}
    \sum_{|n|\in \big[ \frac{\bq_{j+1}}{2^{r+2}},\frac{\bq_{j+1}}{2^{r+1}} \big[} l_n \le 2^{k_j-r}\bq_{j} \exp\left(- \frac{c \, q_j}{j^2} \cdot 2^{(k_j-r)/2}  \right)    
\end{equation}
Now, fix a large-enough $j_0$ such that, for all $j>j_0$, one has 
$$\frac{c \, q_j}{j^2} \ge \sqrt{q_j}\ge 10.$$
Combining this with the inequality $2^{t/2}\ge 1+\frac{t}{10}$, we get that, for $j > j_0$, 
the right-side expression of~\eqref{eq:l-q-r} does not exceed 
\begin{equation}\label{eq:upper-sqrt}
2^{k_j-r}\bq_{j} \exp\left(- \sqrt{\bq_j} \cdot  \left(1+\frac{k_j-r}{10}\right) \right) \le 
2^{k_j-r} e^{-(k_j-r)} \cdot \bq_{j} e^{- \sqrt{\bq_j} }.    
\end{equation}
Now notice that
\[
\bigcup_{r=0}^{k_j} \left[ \frac{\bq_{j+1}}{2^{r+2}},\frac{\bq_{j+1}}{2^{r+1}} \right[
= \left[ \frac{\bq_{j+1}}{2^{k_j+2}},\frac{\bq_{j+1}}{2} \right[ \supset \left[ \frac{\bq_{j}}{2},\frac{\bq_{j+1}}{2} \right[.
\]
Thus, if we sum (\ref{eq:l-q-r}) from $r=0$ to $r=r_j$, using~\eqref{eq:upper-sqrt} we obtain, for $j > j_0$, 
\begin{equation}\label{eq:l-j}
\sum_{|n| \in \big[ \frac{\bq_{j}}{2},\frac{\bq_{j+1}}{2} \big[ } l_n \le 
\sum_{r=0}^{k_j} 2^{k_j-r} e^{-(k_j-r)} \cdot \bq_{j} e^{- \sqrt{\bq_j} }
\le \underbrace{\frac{1}{1-\frac{2}{e}}}_{<4} \cdot \, \bq_{j} e^{- \sqrt{\bq_j} }.    
\end{equation}
The summation over $j>j_0$ thus provides an upper bound over 
\[
|n|\in \bigcup_{j>j_0} \left[ \frac{\bq_{j}}{2},\frac{\bq_{j+1}}{2} \right[ = \Big[ \frac{\bq_{j_0}}{2},+\infty \Big),
\]
thus proving that
\[
\sum_{n\in \Z} l_n \le \sum_{|n|<\bq_{j_0}/2} l_n + 4\sum_{j>j_0} \bq_{j} e^{- \sqrt{\bq_j} } < +\infty.
\]
The finiteness of this sum concludes the proof of Theorem~\ref{t:main}.

\subsection*{Acknowledgements.} The authors are grateful to Étienne Ghys and Maximiliano Escayola for 
helpful discussions and for their interest in this problem, which led to the preparation of the present note.

\end{document}